
\input amstex

\magnification 1200
\loadmsbm
\parindent 0 cm

\define\nl{\bigskip\item{}}
\define\snl{\smallskip\item{}}
\define\inspr #1{\parindent=20pt\bigskip\bf\item{#1}}
\define\iinspr #1{\parindent=27pt\bigskip\bf\item{#1}}
\define\einspr{\parindent=0cm\bigskip}

\input amssym
\input amssym.def

\centerline{\bf A note on Radford's $S^4$ formula}
\bigskip\bigskip
\centerline{\it L.\ Delvaux \rm $^{(1)}$, \it A.\ Van Daele \rm $^{(2)}$ and \it Shuanhong Wang \rm  $^{(3)}$}
\bigskip\bigskip\bigskip
{\bf Abstract} 
\bigskip 
In this note, we show that Radford's formula for the fourth power
 of the antipode can be proven for any regular multiplier Hopf
 algebra with integrals (algebraic quantum groups). This of course not only includes the case of a
 finite-dimensional Hopf algebra but also the case of any
 Hopf algebra with integrals (co-Frobenius Hopf algebras). The
 proof follows in a few lines from well-known formulas in the
 theory of regular multiplier Hopf algebras with integrals.
\snl
 We discuss these formulas and their importance in this theory. We also mention their generalizations, in particular to  the (in a certain sense) more general theory of locally compact  quantum groups. Doing so, and also because the proof of the main result itself is very short, the present note becomes largely of an expository nature.
\nl
\nl
{\it October 2006} (Version 1.2)
\vskip 6 cm
\hrule
\bigskip\parindent 0.5 cm
\item{(1)} Department of Mathematics, University of Hasselt, Agoralaan, B-3590 Diepenbeek (Belgium). E-mail: Lydia.Delvaux\@uhasselt.be
\item{(2)} Department of Mathematics, K.U.\ Leuven, Celestijnenlaan 200B,
B-3001 Heverlee (Belgium). E-mail: Alfons.VanDaele\@wis.kuleuven.be
\item{(3)} Department of Mathematics,  Southeast University,  Nanjing 210096 (China). E-mail: shuanhwang\@seu.edu.cn
\parindent 0 cm

\newpage

\bf 0. Introduction \rm
\nl
Let $H$ be a finite-dimensional Hopf algebra over a field $k$.
 Radford's formula for the fourth power of the antipode says that
$$
S^4(h)=g(\alpha \rightharpoonup h\leftharpoonup \alpha ^{-1})g^{-1}
$$
for any element $h\in H$ where $g$ and $\alpha $ are the so-called
 distinguished group-like elements in $H$ and its dual $H'$
 respectively.
\snl
 The formula was first proven by Larson in a special case [L]
 and later extended by Radford to general finite-dimensional
 Hopf algebras [R]. Since Radford proved the formula in 1976, it
 has been subject to various generalizations. Recently, in an article by Beattie, Bulacu and
 Torrecillas [B-B-T], the formula was proven for any co-Frobenius
 Hopf algebra. See the introduction of that paper for more about the history and recent work on Radford's formula.
\snl
Consider  the finite-dimensional case and observe
 the following. If $\varphi $ is a left integral on $H$, then
 the element $g$ is characterized by the property
 that $\varphi (S(h))=\varphi (hg^{-1})$ for all $h\in H$. Similarly,
 when $\varphi'$ is a left integral on the dual Hopf algebra
 $H'$, the element $\alpha $ is characterized by a similar formula,
 namely $\varphi'(S'(h'))=\varphi '(h'\alpha )$ for all $h'\in H'$
 where $S'$ is the antipode of the dual. We use a slightly
 different terminology than the one common in the
 theory of Hopf algebras. However, doing so, we
 immediately see that Radford's formula can be stated
 for any regular multiplier Hopf algebra with integrals
 $(A, \Delta )$. Indeed, also in this case, we
 have a unique group-like element $\delta $ with the property
 that $\varphi (S(a))=\varphi (a\delta )$ for all $a\in A$. Now $\delta \in M(A)$,
 the multiplier algebra of $A$ (see the preliminaries). Because
 for such a multiplier Hopf algebra, one can construct
 the dual $(\widehat{A}, \widehat{\Delta })$ which is again
 a regular multiplier Hopf algebra with integrals, we also
 have the group-like element $\widehat{\delta }\in M(\widehat{A})$,
 characterized by $\widehat{\varphi }(S(b))=\widehat{\varphi }(b\widehat{\delta })$
 for all $b\in \widehat A$, where now $\widehat{\varphi }$ is a left integral on the dual and
 where $S$ is also used to denote the antipode on the dual. We refer to
 the  literature  and also to the preliminaries for more details about $\delta $
 and $\widehat{\delta }$.
 \snl
 Usually, when a formula is true for finite-dimensional Hopf
 algebras and when it makes sense for regular multiplier Hopf
 algebras with integrals, it is also true in this more general
 setting. And indeed, this turns out to be the case, also
 for Radford's formula. This is the main result discussed in the present note.
 \snl
However, there is more. Not  only do we obtain a substantial
 generalization of the formula, we also show how it follows very
 quickly from other formulas, some of which unpublished, but well-known (for some time)
 to researchers in the field of multiplier Hopf algebras with integrals and
 locally compact quantum groups. Remark that the theory of
 multiplier Hopf algebras was developed first as a by-product
 during the process of the search for a good notion of what
 is now called a locally compact quantum group. This explains the
 fact that the terminology is not completely the same as the one
 used in  Hopf  algebra theory. On the other hand, this disadvantage
 is  compensated by the availability of new techniques.
\snl
Also, in this  note, we will use both the slightly different
 terminology on the one hand, but also other techniques than the
 ones common in Hopf algebra theory. For convenience of the reader,
 we will clearly explain the difference in terminology. And we
 hope the reader will be convinced that these other techniques are
 valuable.
\snl
We have, apart from this introduction, two sections in this paper.
 In Section 1 we will shortly describe the basic notions and
 results in the theory. In Section 2 we will prove the result. We
 will also comment on the various formulas that are used and in
 the mean time, we will refer to various places in literature
 where such formulas can be found, either in special cases, or
 in more general cases.
\snl
Because the proof of the result itself is quite short and because of the focus on these other aspects, the note is mainly expository.
\snl
We refer to Section 1 for basic references, notations and
 terminology.
\nl\nl
\bf Acknowlegdements \rm
\nl
This work was partly supported by the Research Council of the K.U.\ Leuven (through a fellowship for Shuanhong Wang).
\nl\nl

\bf 1. Prelimiaries \rm
\nl
In this preliminary section, we will recall some of the basic notions and results in the theory of regular multiplier Hopf algebras with integrals (sometimes called {\it algebraic quantum groups}). We will be rather short and refer to the literature for more details. However, we will be somewhat more explicit about the differences in terminology between the usual Hopf algebra theory and the theory of multplier Hopf algebras and algebraic quantum groups.
\nl
{\it Multiplier Hopf algebras}
\nl
Let $A$ be an algebra over $\Bbb C$, with or without
 identity, but with a non-degenerate product (i.e.\
 if $a\in A$ and either  $ab=0$ for all $b\in A$ or
 $ba=0$ for all $b\in A$, then $a=0$). If the algebra
 has an identity, of course, the product is automatically
  non-degenerate. The multiplier algebra $M(A)$ of
  $A$ is characterized as the largest algebra with identity,
  containing $A$ as a two-sided ideal with the property
  that still, if $x\in M(A)$ and $xa=0$ for all $a\in A$
  or $ax=0$ for all $a\in A$, then $x=0$. Again, if
  $A$ has an identity then $M(A)=A$.
\snl
Consider the tensor product $A\otimes A$ of $A$ with itself. Then $A\otimes
 A$ is an algebra with a non-degenerate product and we can
 consider the multiplier algebra $M(A\otimes A)$. A coproduct $\Delta $
 on $A$ is a homomorphism from $A$ to $M(A\otimes A)$ satisfying some
 extra conditions (such as coassociativity).
\snl
A {\it multiplier Hopf algebra} is a pair $(A, \Delta )$ of an algebra
 $A$ over $\Bbb C$ with a non-degenerate product and a
 coproduct $\Delta $ such that the two linear maps
 $T_1$, $T_2$ defined from $A\otimes A$ to $M(A\otimes A)$ by
$$\align
 T_1(a\otimes b)&=\Delta (a)(1\otimes b)\\
 T_2(a\otimes b)&=(a\otimes 1)\Delta (b)
\endalign$$
are injective, map into $A\otimes A$ and have range all of $A\otimes A$.
 A multiplier Hopf algebra is called {\it regular} if $(A, \Delta ')$
 is still a multiplier Hopf algebra with $\Delta '$ obtained
 from $\Delta $ by composing it with the flip.
\snl
For any multiplier Hopf algebra $(A, \Delta )$ there is
 a counit $\varepsilon$ and an antipode $S$ satisfying
 (and characterized by) properties, very similar
  as in the case of Hopf algebras. Regularity
  is equivalent with the antipode being a bijective map
  from $A$ to itself.
\snl
Any Hopf algebra is a multiplier Hopf algebra. Conversely,
 if $(A, \Delta )$ is a multiplier Hopf algebra and if $A$ has an
 identity, then it is a Hopf algebra.
\snl
The motivating example comes from a group $G$. Let
 $A$ be the algebra $K(G)$ of complex functions on $G$
  with finite support and pointwise product. Identify
  $A\otimes A$ with $K(G\times G)$ and $M(A\otimes A)$ with the algebra
  $C(G\times G)$ of all complex functions on $G\times G$. Then
  $A$ becomes a (regular) multiplier Hopf algebra if the coproduct
  $\Delta $ is defined by $\Delta (f)(p, q)=f(pq)$ for $p, q\in G$
 and $f\in K(G)$.
\snl
The theory of multiplier Hopf algebras has been developed for algebras over the field $\Bbb C$ (because of an analytical background) but this is not essential and any other field would do as well.
\snl
We refer to [VD1] for details about the theory of multiplier
 Hopf algebras. For more examples, see e.g.\ [VD2], [VD-Z], and for some constructions, see e.g.\ [Dr-VD], [D-VD3] and [D-VD-W1].
\nl
{\it Multiplier Hopf algebras with integrals}
\nl
Assume in what follows that $(A, \Delta )$ is a regular multiplier
 Hopf algebra. A linear functional $\varphi$ on $A$ is called
 left invariant if  $(\iota\otimes \varphi)\Delta (a)=\varphi(a)1$ in $M(A)$
 for all $a\in A$. We use $\iota$ to denote the identity map
 and the formula is given a meaning by multiplying from left
 or right with any element in $A$. A non-zero left invariant
 functional $\varphi$ is called a {\it left integral on} $A$. Similarly,
  a right integral is defined. The motivating example, explained
  before, justifies the terminology. Indeed, if $A=K(G)$ for a
  group $G$ and $\varphi$ is defined by $\varphi(f)=\sum_{p\in G}f(p)$,
  then $\varphi$ is a left (and a right) integral.
\snl
In general, one can show that left and right integrals are unique
 (up to a scalar) if they exist. If a left integral $\varphi$ exists also a right integral
 exists, namely $\varphi\circ S$.  Left and right
  integrals may be different but they are related. For a left integral $\varphi$, there is a
  unique group-like element $\delta$ in $M(A)$ such that $\varphi(S(a))=\varphi(a\delta)$
  for all $a\in A$. It is (the inverse of) the
  so-called distinguished group-like element, known in Hopf
  algebra theory. In the theory of multiplier Hopf algebras, it is
  called the {\it modular element} because it relates the left and the
  right integral just as in the theory of locally compact groups,  the modular function relates the left
 and the right Haar integrals.
\snl
Let $\varphi$ be a left integral and $\psi$ a right integral. 
There are automorphisms $\sigma $ and $\sigma '$ of $A$ satisfying
$$\align \varphi(ab)&=\varphi(b\sigma (a))\\
\psi (ab)&=\psi (b\sigma '(a))
\endalign$$
for all $a, b\in A$. It is known that the integrals are
 faithful (i.e.\ if $a\in A$ and $\varphi(ab)=0$ for all $b\in A$
 or $\varphi(ba)=0$ for all $b\in A$, then $a=0$; and similarly
 for $\psi $). So these automorphisms are unique and characterized
 by these formulas. The property is called {\it the weak K.M.S.\
 property} and the automorphisms are called the {\it modular
 automorphisms}. This terminology comes from the theory of
 operator algebras and finds its origin in physics. In Hopf algebra theory, these automorphisms are known as the Nakayama automorphisms.
\snl
Finally, there is a scalar $\tau \in \Bbb C$ satisfying
 $\varphi(S^2(a))=\tau \varphi(a)$ for all $a\in A$. It exists because
 $\varphi\circ S^2$ is also a left integral and therefore a multiple
 of $\varphi$. For obvious reasons, it is called the {\it scaling constant}.
\snl
There are many formulas relating these objects. We will state them when we use them in the next section. 
For details about regular multiplier Hopf algebras with integrals we refer to [VD2], see also
 [VD-Z].  The example $K(G)$ is too trivial to illustrate these
 various objects. See [VD2] and [VD-Z] for examples that do illustrate
 various features.
\nl
{\it Duality for multiplier Hopf algebras with integrals}
\nl
In what follows, $(A, \Delta )$ is a regular multiplier Hopf algebra with integrals and
 $\varphi$ and $\psi $ are a left and a right integral respectively.
\snl
Define $\widehat {A}$ as the space of linear
 functionals on $A$ of the form $\varphi(\cdot\, a)$
 where $a\in A$. Then $\widehat{A}$
 can be made into an algebra with a product dual to the coproduct
 $\Delta $ on $A$. This product is non-degenerate.
 The product on $A$ can also be dualized and yields a
 coproduct $\widehat {\Delta }$ on $\widehat {A}$. It is shown
 that the pair $(\widehat {A}, \widehat {\Delta })$
 becomes a regular multiplier Hopf algebra. It is
 called the {\it dual of} $(A, \Delta )$.
\snl
There exist integrals on $(\widehat {A}, \widehat {\Delta })$. A right
 integral $\widehat {\psi }$ on $\widehat {A}$ is defined by
 $\widehat {\psi }(\omega)=\varepsilon(a)$ if $\omega=\varphi(\cdot \, a)$ and $a\in A$.
 The various objects associated with $(\widehat {A}, \widehat {\Delta })$
 are denoted as for $(A, \Delta )$ but with a {\it hat}. So
 we have $\widehat {\varphi}$, $\widehat {\psi }$, $\widehat {\delta}$, $\widehat {\sigma }$,
 $\widehat {\sigma }'$. However, we use $\varepsilon$ and $S$ also for
  the counit and antipode on the dual.
\snl
Repeating the procedure, i.e. if we take the dual of $(\widehat {A},
 \widehat {\Delta })$, we find the original pair $(A, \Delta )$. The result
 is not very deep. It essentially follows from the formula
 $$\widehat \psi(\omega'\omega)=\omega'(S^{-1}(a))$$
for all $a\in A$ and $\omega'\in \widehat A$, where $\omega=\varphi(\cdot \, a)$. This formula in turn follows quite
 easily from the definitions (see e.g.\ Theorem 4.12 in [VD2]). The result is refered to as {\it biduality}.
\snl
Duality for regular multiplier Hopf algebra with integrals
 generalizes duality for finite-dimensional Hopf algebras.
 The results are quite similar but of course, the class of objects
 is much bigger. It can be applied to Hopf algebras with integrals
 (the co-Frobenius type) but also to the duals. In the theory of
  multiplier Hopf algebras, we call the first multiplier Hopf algebras of compact type
   and the second multiplier Hopf algebras of discrete type.
   The terminology has its origin in analysis. Think of  Pontryagin
   duality for abelian locally compact groups. The dual of a
   compact group is a discrete group. And to any compact
   group $G$ is associated  the Hopf algebra of polynomial functions
   on $G$. For a discrete group, we have the multiplier Hopf
   algebra $K(G)$ as introduced before. It is indeed of discrete
 type.

\snl
This is perhaps a good occasion to say something more about terminology related with integrals. For us, an integral is always a linear functional {\it on} the algebra. This is because the algebra is considered as an 'algebra of functions' on a 'quantum space'. And in classical analysis, an integral is a linear functional on a space of functions. What is often called an integral {\it in} a Hopf algebra, we will call a {\it co-integral} because it should be considered as an integral {\it on} the dual Hopf algebra. We use this terminology consistently in the theory of algebraic quantum groups. For an algebraic quantum group $(A,\Delta)$, an integral is a linear functional on $A$. If the integral on the dual $\widehat A$ turns out to be an element in $A$, then we call it a co-integral in $A$. This is precisely the case when $A$ is of discrete type.   
\snl
For more details about duality, we again refer to [VD2]. 
\snl
Also here, there are many
 formulas relating the various objects for $(A, \Delta )$ with these
 of $(\widehat {A}, \widehat {\Delta })$. Only a few of them are given in [VD2]. More formulas are found in [K], but it should be mentioned that these are only proven in the case of a multiplier Hopf $^*$-algebra with positive integrals. However, in general, when a formula makes sense in this special case, it turns out to be true, also in the general case.
We will recall (some of) these formulas (and prove them if necessary) where we use them in the next section. 
\nl
\it Pairing and actions \rm
\nl
We will use the expression $\langle a,b \rangle$ to denote the value of a functional $b\in \widehat A$ in the point $a\in A$. Doing so, we get a non-degenerate pairing of multiplier Hopf algebras in the sense of [Dr-VD]. This pairing gives rise to the natural four actions. We have a left and a right action of $A$ on $\widehat A$, given by the formulas
$$\align \langle a'a,b \rangle &= \langle a', a\rightharpoonup b\rangle \\
\langle aa', b\rangle&=\langle a', b\leftharpoonup a \rangle
\endalign$$
for any $a,a'\in A$ and $b\in \widehat A$. Similarly the left and right actions of $\widehat A$ on $A$ are given by 
$$\align \langle a,b'b\rangle&=\langle b\rightharpoonup a,b'\rangle \\
\langle a,bb' \rangle &=\langle a\leftharpoonup b, b'\rangle
\endalign$$
for any $a\in A$ and $b,b'\in \widehat A$. It is not completely obvious that this can be done, but we refer to [Dr-VD] for a precise treatment. See also [Dr-VD-Z] for more information about actions of multiplier Hopf algebras in general.
\nl\nl

\bf 2. The main result and comments \rm
\nl
In this section, the aim is twofold. On the one hand, we will
 present a proof of the main result, Radford's $S^4$ formula for
 regular multiplier Hopf algebras with integrals (algebraic quantum
 groups). On the other hand, we will also comment on the various
 formulas, related with and necessary for the proof of  Radford's
 formula. We will make links with other results in the field and compare with 
results in the (analytic) theory of locally compact quantum groups.
\snl
 We consider a regular multiplier Hopf algebra
 $(A, \Delta )$ with integrals  and its dual $(\widehat   {A}, \widehat   {\Delta })$ as
 explained in the previous section.
\snl
We have the different objects associated with $(A, \Delta )$ and $(\widehat  
 {A}, \widehat   {\Delta })$ and many relations among them. In particular, we
 have a left and a right integral $\varphi   $ and $\psi $ on $A$ and we
 have a left and a right integral $\widehat   {\varphi   }$ and $\widehat   {\psi }$ on $\widehat A$.
 We fix $\varphi   $ and we normalize the three others in the following
 way. We let $\psi =\varphi  \circ S$ and $\widehat   {\psi }=\widehat   {\varphi  }\circ
 S$. We define $\widehat   {\psi }$ by $\widehat   {\psi }(\omega)=\varepsilon(a)$ when
  $\omega=\varphi  (\cdot \, a)$ and $a\in A$ as before. One can show that,
  with these assumptions, we have $\widehat   {\varphi  }(\omega)=\varepsilon(a)$ when
  $\omega= \psi (a\,\cdot \,)$ and $a\in A$ (see Proposition 4.8 in [VD2]). It should be mentioned however that for our purposes, this normalization is not so important.
\snl
For the modular automorphism $\sigma$ and $\sigma'$ we have the two
 basic formulas
$$\align
\Delta (\sigma(a))&=(S^2\otimes\sigma)\Delta (a)\\
\Delta (\sigma'(a))&=(\sigma'\otimes S^{-2})\Delta (a)
\endalign$$
for all $a\in A$. They are found already in [VD2, Proposition 3.14]. However, there
 is also a third formula of this type. It says that, not only 
$\Delta (S^2(a))=(S^2\otimes S^2)\Delta (a)$, but also 
$$
\Delta (S^2(a))=(\sigma\otimes{\sigma'}^{-1})\Delta (a)
$$
for all $a\in A$. This formula was first proven in [K-VD1, Lemma 3.10], in the
 case of a multiplier Hopf $^*$-algebra with positive integrals, but later simpler and more direct proofs have been
 given, valid in the general case (see the appendix in  [L-VD] and also Proposition 2.7 in [D-VD2]). The
 automorphism $\sigma, \sigma'$ and $S^2$ all mutually commute and
 there are two basic relations between $\sigma$ and $\sigma'$. We have
 $S(\sigma'(a))=\sigma^{-1}(S(a))$ and also $\delta \sigma(a)=\sigma'(a)\delta$
 for all $a\in A$. These properties are already found in the
 original paper [VD2].
\snl
Now, we also have the data for the dual and of course, all
 the above results can also be formulated for the dual objects. In
 the original paper however, very few relations between the
 objects for $(A, \Delta )$ and the ones for the dual $(\widehat   {A}, \widehat   {\Delta })$
 are found. Many such relations (and much more) are found in [K].
 Unfortunately,  in that paper, it is assumed that $(A, \Delta )$ is a
 multiplier Hopf $^*$-algebra (i.e.\ $A$ is a $^*$-algebra and
 $\Delta $ a $^*$-homomorphism) and the left integral is assumed to
 be positive (i.e.\ $\varphi  (a^*a)\geq 0$ for all $a$). Remember
 that in that case, the right integral $\psi $ is automatically
 positive. Therefore, as we mentioned already, the results of [K] can
 not be used in the general case we consider here. On the other
 hand, these results are a good source of inspiration and often they are also true in the more general case. 
\nl
The first of such formulas we like to consider and discuss are the
 following.

\inspr{2.1} Proposition \rm
 For all $a\in A$ we have
$$\align
&\langle \, a, \, \widehat   {\delta } \, \rangle = \varepsilon(\sigma^{-1}(a))=\varepsilon({\sigma'}^{-1}(a))\\
&\langle \, a, \widehat   {\delta }^{-1}\, \rangle =\varepsilon(\sigma(a))=\varepsilon(\sigma'(a)).
\endalign$$
\einspr
Before we give a proof, we need to make a comment on the left hand
 side of these formulas. We need to observe that the original
 pairing of $A$ with $\widehat{A}$ can be extended in a natural
 way to a pairing between $A$ and $M(\widehat{A})$. If e.g.\ $m$ is multiplier in $M(\widehat A)$, then formally we have
$$\langle a,mb \rangle = \langle b \rightharpoonup a, m \rangle$$
for any $a\in A$ and $b\in \widehat A$. This formula can be used to extend the pairing because the action of $\widehat A$ on $A$ is unital, i.e.\ for all $a\in A$ there is a $b\in \widehat A$ satisfying $b \rightharpoonup a=a$ (see e.g.\ [Dr-VD]. For a precise treatment of this extended pairing, see e.g.\ [D].
\snl
Now, we give the proof of Proposition 2.1.
\inspr{} Proof: \rm
Take $c\in A$ and let $b=\varphi  (\cdot \, c)$. We will use the Sweedler notation (as justified in e.g.\ [Dr-VD] and [Dr-VD-Z]). We write, for any $a, a'\in
 A$,
$$
\varphi  (a' c\sigma(a))=\varphi  (aa'c)=\langle aa', b  \rangle =\langle  a, b _{(1)}\rangle \langle  a', b _{(2)}\, \rangle 
$$
so that
$$
\varphi  (\cdot \, c\sigma(a))=\langle  a, b _{(1)}\, \rangle  b _{(2)}= b \leftharpoonup a.
$$
From the definition of $\widehat   {\psi }$ and the formula $(\iota\otimes\widehat\psi)\Delta(b)=\widehat\psi(b) \widehat \delta^{-1}$, we get
$$
\varepsilon(c\sigma(a))=\widehat   {\psi }(b_{(2)})\langle a, b _{(1)} \rangle =\langle \, a,
\widehat   {\delta}^{-1} \rangle \widehat   {\psi }(b)=\langle  a, \widehat   {\delta }^{-1} \rangle \varepsilon(c)
$$
so that $\varepsilon(\sigma(a))=\langle  a, \widehat   {\delta}^{-1}\, \rangle $.
\snl
This proves one of the formulas. Using the two relations between $\sigma$ and $\sigma'$ formulated earlier, and using that $\varepsilon\circ S=\varepsilon$, $\varepsilon(\delta)=1$ and $S(\widehat   {\delta })=\widehat   {\delta }^{-1}$,
  we easily obtain the other three. \hfill $\blacksquare$
\einspr

From these formulas, we see that $\widehat \delta$ has to be group-like because $\varepsilon$ and $\sigma$ are algebra maps.
\snl
In the original paper [VD2], it was remarked already
 that $\varepsilon(\sigma(a))=\varepsilon(\sigma'(a))$ for all $a$ but it was also stated {\it 
 There seems to be no obvious relation between $\sigma$ and $\varepsilon$} (see the end of Section 3 in [VD2]). However, these formulas are
 essentially obtained in [K], but as said, only for
 multiplier Hopf $^*$-algebras with positive integrals. Indeed, Proposition 7.10 in [K] says that $\widehat \delta^{iz}=\varepsilon\sigma_{-z}$ and this will give the results in our Proposition 2.1 above when we let $z=i$ and $z=-i$.
It was known however for some time that these formulas also were true in the general case. They are found e.g.\ in an unpublished preprint [K-VD2]. In [VD-W], they were obtained in a special case. Recently, they have been included in [D-VD2] for the more general case of algebraic quantum hypergroups.
\nl
But, as we saw, the formulas in Proposition 2.1 are easy to obtain. Nevertheless, it are precisely these formulas that will give us very quicly, Radford's formula.
\snl
First we need to prove some other formulas. These are equivalent with the formulas in the previous proposition in the sense 
that it is easy to obtain one set from the other.

\inspr{2.2} Proposition \rm
 For all $a\in A$ and $b\in \widehat   {A}$, we
 have
$$\align
&\langle \sigma(a),b  \rangle         =\langle a,   S^2(b)\widehat\delta^{-1} \rangle\\
&\langle\sigma^{-1}(a),b \rangle   =\langle a, S^{-2}(b)\widehat\delta \rangle\\
&\langle\sigma'(a),b \rangle           =\langle a, \widehat\delta^{-1}S^{-2}(b) \rangle\\
&\langle{\sigma'}^{-1}(a),b)\rangle=\langle a, \widehat\delta S^2(b) \rangle.
\endalign$$

\snl\bf Proof: \rm Start with the formula $\Delta (\sigma(a))=(S^2\otimes\sigma)\Delta
 (a)$. If we apply $\iota\otimes\varepsilon$ and use the previous proposition, we get
 $$
 \sigma(a)=S^2(a_{(1)})\varepsilon(\sigma(a_{(2)}))
  =\langle \, a_{(2)}, \widehat   {\delta}^{-1}\, \rangle S^2(a_{(1)}).
 $$
Pairing with $b$ gives
$$
 \langle \, \sigma(a), b\, \rangle=\langle \, a, S^2(b)\widehat   {\delta}^{-1}\, \rangle.
$$
This gives the first formula. The three other ones are obtained in a similar fashion.\hfill $\blacksquare$
\einspr

If in the formulas above, we let $b=1$, we see that we recover the formulas in the previous proposition. One has to be a little bit careful, but this is done by using the 'covering technique', needed to work properly with Sweedler's notation in the context of multiplier Hopf algebras (cf.\ [Dr-VD] and [Dr-VD-Z]). 
\snl
The formulas in Proposition 2.1 and in Proposition 2.2 have obvious dual forms (see e.g.\ [K] and [D-VD2]). The forms we have given are the ones that are needed to prove Radford's formula.
\snl
The formulas are very useful. We see e.g.\ that, if $\varphi$ is a trace, that is when $\sigma$ is trivial, then it follows from Proposition 2.1 that $\widehat \delta=1$ (meaning that left and right integrals on the dual are the same, i.e.\ the dual is {\it unimodular}). The first formula in Proposition 2.2 will give us in this case that $S^2=\iota$. In particular, when $A$ is abelian, we recover the result that the square of the antipode is trivial. By duality, the same is true when $A$ is coabelian. Another conclusion from these formulas is that, if left and right integrals on $A$ coincide (i.e.\ that $A$ is unimodular and $\delta=1$), then $\widehat\sigma$ and $S^2$ coincide on $\widehat A$. 
\snl
Again, these formulas can be found already in [K]. Consider Proposition 7.3 in that paper where it is stated that 
$$\langle a, \widehat\sigma_z(b)\rangle= \langle \tau_z(a)\delta^{-iz},b\rangle$$
whenever $a\in A$ and $b\in \widehat A$. If we take $z=-i$ we get
$$\langle a, \widehat\sigma(b)\rangle= \langle S^2(a)\delta^{-1},b\rangle$$
because $\widehat\sigma_{-i}$ and $\tau_{-i}$ in [K] are here $\widehat \sigma$ and $S^2$ respectively. So, this gives the first formula, in dual form, of Proposition 2.2 above. The formulas were also discovered in the unpublished preprint [K-VD2] and finally proven in the more general setting of algebraic quantum hypergroups in [D-VD2]. Related results were also obtained in the unpublished paper [VD-W], but only for multiplier Hopf algebras of discrete type.
\snl
In the analytic setting of locally compact quantum groups, we have equivalent formulas. Take e.g.\ the first formula in Theorem 4.17 of [VD4]. It comes in the form
$$\nabla^{it}=(\widehat J \widehat \delta^{it}\widehat J)P^{it}$$
for all $t\in \Bbb R$. This is an equality of unitary operators on a Hilbert space. With $t=i$ and $t=-i$, they are essentially the same as the first two formulas in Proposition 2.2 here. It is not so easy to see this, and one should notice that the operator $\nabla$ is related with the automorphism $\sigma$ whereas the operator $P$ is related with the square of the antipode. These formulas are already found in the original papers by Kustermans and Vaes on locally compact quantum groups (see e.g.\ [K-V3]). We do not have counterparts of the formulas of Proposition 2.1 because the counit is, from an analytical point of view, a difficult and not very usefull object in the theory of locally compact quantum groups.
\nl
We now arrive at the main result of this note: Radford's formula for the fourth power of the antipode.

\inspr{2.3} Theorem \rm Let $(A,\Delta)$ be a regular multiplier Hopf algebra with integrals (an algebraic quantum group). If $\delta$ and $\widehat \delta$ denote the modular elements in $M(A)$ and $M(\widehat A)$ respectively, then
$$S^4(a)=\delta^{-1}(\widehat \delta\rightharpoonup a\leftharpoonup \widehat\delta^{-1})\delta$$
for all $a\in A$.

\snl\bf Proof: \rm
From the previous proposition, we see that
$$\align \sigma(x)&=\widehat\delta^{-1}\rightharpoonup S^2(x)\\
              \sigma'(x)&=S^{-2}(x)\leftharpoonup \widehat\delta^{-1}
\endalign$$
for all $x\in A$. If we use that $\delta\sigma(x)=\sigma'(x)\delta$ and then replace $x$ by $S^2(a)$, we find
$$(a\leftharpoonup \widehat\delta^{-1})\delta=\delta (\widehat\delta^{-1}\rightharpoonup S^4(a))$$
for all $a$. Now, let $\widehat \delta$ act on this equation from the left. Because $\Delta(\delta)=\delta\otimes \delta$ and $\varepsilon(\sigma^{-1}(\delta))=\tau$ (the scaling constant), this action commutes, up to $\tau$, with both multiplication from left or right by $\delta$. Then we get  Radford's formula.\hfill $\blacksquare$
\einspr

The first thing we would like to remark is that Radford's formula is {\it self-dual}, in the following sense. If we pair the formula with any element $b\in \widehat A$, we can move $S^4$ to the other side (because $\langle S(a),b\rangle=\langle a,S(b)\rangle$), we can move $\delta^{-1}$ and $\delta$ to the other side to get $\delta\rightharpoonup b\leftharpoonup\delta^{-1}$ and finally, we can move also $\widehat \delta$ and $\widehat\delta^{-1}$ and obtain $\widehat \delta^{-1}(\delta\rightharpoonup b\leftharpoonup \widehat\delta^{-1})\widehat\delta$. So we again get Radford's formula, now for the dual $\widehat A$.
\snl
In [B-B-T] the formula is obtained for a Hopf algebra with integrals. In [VD-W], the formula was discovered already for any regular multiplier Hopf algebra with integrals, but with $\alpha^{-1}$ and $\alpha$ in stead of $\widehat\delta$ and $\widehat\delta^{-1}$ where $\alpha=\varepsilon\sigma$. In that paper, it was noticed that $\alpha=\widehat\delta^{-1}$, but only for multiplier Hopf algebras of discrete type, i.e.\ the dual type of Hopf algebras with integrals (as in [B-B-T]). In [D-VD1], the formula was obtained for general multiplier Hopf algebras with integrals. This was mentioned at a talk of the second author at the AMS meeting in Bowling Green in 2005 [VD3]. In [D-VD2] however, Radford's formula is obtained for any algebraic quantum hypergroup, including the case of any algebraic quantum group.
\snl
In the case of a multiplier Hopf $^*$-algebra with positive integrals, we have analytic forms of all the formulas in Proposition 2.1 and Proposition 2.2 thanks to the work of Kustermans [K]. This means that there is also an analytic form of Radford's formula in this case. We use this to construct a ribbon element in [D-VD-W2]. In fact, in [VD4], we proved an analytic version of Radford's formula for any locally compact quantum group. It comes under the form
$$P^{-2it}=\delta^{it}(J\delta^{it}J)\widehat\delta^{it}(\widehat J\widehat\delta^{it}\widehat J)$$
for all $t\in \Bbb R$ (cf.\ Theorem 4.20 in [VD4]). Here again, we have an equation with unitary operators on a Hilbert space. It is not completely obvious, but with $t=i$, we recover a formula (of undbounded operators) that is essentially Radford's formula. The above formula is not found in the original papers by Kustermans and Vaes ([K-V1] and [K-V2]), but it follows quickly from formulas found in [K-V3] (see an earlier remark).
\snl
Finally, we also want to refer to the work of C.\ Voigt. In [V1] he introduced the notion of a bornological quantum group. Basically, it is a bornological vector space, endowed with a suitable product and coproduct. Regular multiplier Hopf algebras with integrals are a special case. In fact, the theory of algebraic quantum groups served as a model for the development of these bornological quantum groups. In Section 4 of [V2], a proof is given of Radford's formula for a bornological quantum group. The formulation is just a for algebraic quantum groups (as in Theorem 2.3 of this note) and also the proof is very similar. As bornological quantum groups are generalizing the algebraic quantum groups, the result of C.\ Voigt is more general than our Theorem 2.3. 
\nl\nl

\bf References \rm
\nl
\parindent 1 cm

\item{[B-B-T]} M.\ Beattie, D.\ Bulacu \& B.\ Torrecillas: {\it Radford's $S^4$ formula for co-Frobenius Hopf algebras}. J.\ Algebra (2006), to appear. 
\item{[D]} L.\ Delvaux: {\it The size of the intrinsic group of a multiplier Hopf algebra}. Commun.\ Alg.\ 31 (2003), 1499--1514.
\item{[D-VD1]} L.\ Delvaux \& A.\ Van Daele: {\it Traces on (group-cograded) multiplier Hopf algebras}. Preprint University of Hasselt and K.U.\ Leuven (2005).
\item{[D-VD2]} L.\ Delvaux \& A. Van Daele: {\it Algebraic quantum hypergroups}. Preprint University of Hasselt and K.U.\ Leuven (2006). math.RA/0606466
\item{[D-VD3]} L.\ Delvaux \& A.\ Van Daele: {\it The Drinfel'd double of multiplier Hopf algebras}. J.\ Alg.\ 272 (2004), 273--291.
\item{[D-VD-W1]} L.\ Delvaux, A. Van Daele \& Shuanhong Wang: {\it Bicrossed product of multiplier Hopf algebras}. Preprint University of Hasselt, K.U.\ Leuven and Southeast University Nanjing (2005).
\item{[D-VD-W2]} L.\ Delvaux, A. Van Daele \& Shuanhong Wang. {\it Quasi-triangular and ribbon multiplier Hopf algebras} Preprint University of Hasselt, K.U.\ Leuven and Southeast University Nanjing. In preparation.
\item{[Dr-VD]} B.\ Drabant \& A. Van Daele: {\it Pairing and Quantum double of multiplier Hopf algebras}.  Algebras and Representation Theory 4 (2001), 109-132.
\item{[Dr-VD-Z]} B.\ Drabant, A.\ Van Daele \& Y.\ Zhang: {\it Actions of multiplier Hopf algebras}. Commun.\ Alg.\ 27 (1999), 4117--4172.
\item{[K]} J.\ Kustermans: {\it The analytic structure of algebraic quantum groups}. J.\ Algebra 259 (2003), 415--450.
\item{[K-V1]} J.\ Kustermans \& S.\ Vaes: {\it A simple definition for locally compact quantum groups}. C.R.\ Acad.\ Sci.\ Paris S\'er I 328 (1999), 871--876.
\item{[K-V2]} J.\ Kustermans \& S.\ Vaes: {\it Locally compact quantum groups}. Ann.\ Sci.\ \'Ecole Norm.\ Sup. (4) (33) (2000), 837--934.
\item{[K-V3]} J.\ Kustermans \& S.\ Vaes: {\it Locally compact quantum groups in the von Neumann algebra setting}. Math.\ Scand.\ 92 (2003), 68--92. 
\item{[K-VD1]} J.\ Kustermans \& A.\ Van Daele: {\it C$^*$-algebraic quantum groups arising from algebraic quantum groups}. Int.\ J.\ Math.\ 8 (1997), 1067--1139.
\item{[K-VD2]} J.\ Kustermans \& A.\ Van Daele: {\it The Heisenberg commutation relations for an algebraic quantum group}. Preprint K.U.\ Leuven (2002), unpublished.
\item{[L]} R.G.\ Larson: {\it Characters of Hopf algebras}. J.\ Algebra 17 (1971), 352--368.
\item{[L-VD]} M.B.\ Landstad \& A.\ Van Daele: {\it Compact and discrete subgroups of algebraic quantum groups}. Preprint University of Trondheim and K.U.\ Leuven (2006).
\item{[R]} D.\ Radford: {\it The order of the antipode of any finite-dimensional Hopf algebra is finite}. Amer.\ J.\ Math.\ 98 (1976), 333--355. 
\item{[VD1]} A.\ Van Daele: {\it Multiplier Hopf algebras}. Trans.\ Am.\ Math.\ Soc.\ 342 (1994), 917--932.
\item{[VD2]} A.\ Van Daele:  {\it An algebraic framework for group duality}.  Adv.\ Math.\ 140 (1998), 323--366.
\item{[VD3]} A.\ Van Daele: {\it Multiplier Hopf algebras with integrals}. Talk at the AMS meeting in Bowling Green, Kentucky (2005).
\item{[VD4]} A.\ Van Daele: {\it Locally compact quantum groups. A von Neumann algebra approach}. Preprint K.U.\ Leuven (2006), math.OA/0602212
\item{[VD-W]} A.\ Van Daele \& Shuanhong Wang: {\it Modular multipliers of regular multiplier Hopf algebras of discrete type}. Preprint K.U.\ Leuven (2004), unpublished.
\item{[VD-Z]} A.\ Van Daele \& Y.\ Zhang: {\it A survey on multiplier Hopf algebras}. In 'Hopf algebras and Quantum Groups', eds. S.\ Caenepeel \& F.\ Van Oyestayen, Dekker, New York (1998), pp. 259--309.
\item{[V1]} C.\ Voigt: {\it Bornological quantum groups}. Preprint University of Copenhagen (2005), math.QA/0511195. 
\item{[V2]} C.\ Voigt: {\it Equivariant cyclic homology for quantum groups}. Preprint University of Copenhagen (2006), math.KT/0601725.
\end